# DADS Under Unknown Input Coefficients

**Iasson Karafyllis[*] and Miroslav Krstic[**]**

[*]Dept. of Mathematics, National Technical University of Athens, Zografou Campus, 15780, Athens, Greece,
email: iasonkar@central.ntua.gr; iasonkaraf@gmail.com

[**]Dept. of Mechanical and Aerospace Eng., University of California, San Diego, La Jolla, CA 92093-0411, U.S.A., email: krstic@ucsd.edu

## Abstract

This short note shows that the Deadzone-Adapted Disturbance Suppression (DADS) adaptive control scheme is applicable to systems with unknown input coefficients. We study time-invariant, control-affine systems that satisfy the matching condition for which no bounds for the disturbance and the unknown parameters are known. The input coefficients can be time-varying as well as the unknown parameters. The only thing assumed for the input coefficients is their sign. The adaptive control design is Lyapunov-based and can be accomplished for every system for which a smooth globally stabilizing feedback exists when the disturbances are absent and all unknown parameters are known. The design is given by simple, explicit formulas. The proposed controllers guarantee an attenuation of the plant state to an assignable small level, despite unknown bounds on the parameters and disturbance, without a drift of the gain, state, and input.

**Keywords:** Robust Adaptive Control, Feedback Stabilization, Matching Condition, Leakage.

## 1. Introduction

Robustness in adaptive control is an important issue that has attracted the attention of many control theorists. The literature offers numerous approaches for time-invariant nonlinear control systems for which no persistence of excitation is assumed: leakage (see [4, 5, 23, 26]), nonlinear damping (see [17, 18, 10, 11]), projection methodologies (see [4] and Appendix E in [17]), dynamic (high) gains or gain adjustment (see [6, 16, 22]), supervision for direct adaptive schemes (see [3]), and deadzone in the update law (introduced in the paper [24] and well explained in the book [4]). Every approach has its own strong points while some of the approaches require special assumptions (e.g. knowledge of bounds for the disturbances and/or the unknown parameters). However, none of them achieve convergence of the plant state to an assignable small neighborhood of the origin, despite unknown bounds on the parameters and disturbance, without a drift of the gain, state, and input. Indeed, among the robustification tools for the update law, only sigma-modification (see [4, 5, 23, 26]) is free of a requirement of a known bound on the parameter. Projection requires and bound and the deadzone alone is known to provide only local robustness (local in the initial condition and the disturbance size).

Recently, a novel adaptive control scheme for time-invariant systems was proposed in [12, 13, 14, 15]: the Deadzone-Adapted Disturbance Suppression (DADS) control scheme. Inspired by the idea of dynamic gain used in [6, 25, 16, 21, 22], DADS provides a simple, direct, adaptive control scheme that combines three elements: (a) nonlinear damping (as in [17, 11]), (b) single-gain adjustment (the dynamic feedback has only one state), and (c) deadzone in the update law. The proposed adaptive scheme is direct and robustness is not sought by applying advanced identification tools or delays (as in [8, 9]); no identification is performed (which also explains the simplicity of the proposed scheme). The DADS controller achieves an attenuation of the plant state to an assignable small level, despite the presence of disturbances and unknown parameters of arbitrary and unknown bounds. Moreover, the DADS controller prevents gain and state drift regardless of the size of the disturbance and unknown parameter.

In the present short note, we focus on a special class of systems: time-invariant systems that satisfy the so-called matching condition (see [17] for the definition of the matching condition) for which no bounds for the disturbance and the unknown parameters are known. We show that the DADS control scheme is applicable to systems with unknown input coefficients. The input coefficients can be time-varying as well as the unknown parameters. The only thing assumed for the input coefficients is their sign. The DADS design is Lyapunov-based and can be accomplished for every control-affine system for which a smooth globally stabilizing feedback exists when the disturbances are absent and all unknown parameters are known. The design is given by simple, explicit formulas (see Section 2 below).

The structure of the paper is as follows. We start with a subsection that provides all the stability notions used in the paper. All the main results are stated and discussed in Section 2. Section 3 is devoted to the presentation of a simple pedagogical example: the uncertain double integrator. Section 4 of the paper contains the proof of the main result.

**Notation and Basic Notions.** Throughout this paper, we adopt the following notation.

* $\mathbb{R}_+ := [0, +\infty)$. For a vector $x \in \mathbb{R}^n$, $|x|$ denotes its Euclidean norm and $x'$ denotes its transpose. We use the notation $x^+$ for the positive part of the real number $x \in \mathbb{R}$, i.e., $x^+ = \max(x, 0)$.

* Let $D \subseteq \mathbb{R}^n$ be an open set and let $S \subseteq \mathbb{R}^n$ be a set that satisfies $D \subseteq S \subseteq cl(D)$, where $cl(D)$ is the closure of $D$. By $C^0(S; \Omega)$, we denote the class of continuous functions on $S$, which take values in $\Omega \subseteq \mathbb{R}^m$. By $C^k(S; \Omega)$, where $k \geq 1$ is an integer, we denote the class of functions on $S \subseteq \mathbb{R}^n$, which take values in $\Omega \subseteq \mathbb{R}^m$ and have continuous derivatives of order $k$. In other words, the functions of class $C^k(S; \Omega)$ are the functions which have continuous derivatives of order $k$ in $D = \text{int}(S)$ that can be continued continuously to all points in $\partial D \cap S$. When $\Omega = \mathbb{R}$ then we write $C^0(S)$ or $C^k(S)$. A function $f \in \bigcap_{k=0}^{\infty} C^k(S; \Omega)$ is called a smooth function.

* Let $D \subseteq \mathbb{R}^p$ be a non-empty set. By $L^\infty(\mathbb{R}_+; D)$ we denote the class of essentially bounded, Lebesgue measurable functions $d: \mathbb{R}_+ \to D$. When $D = \mathbb{R}^p$ then we simply write $L^\infty(\mathbb{R}_+)$. For $d \in L^\infty(\mathbb{R}_+; D)$ we define $\|d\|_\infty = \sup_{t \geq 0}(|d(t)|)$, where $\sup_{t \geq 0}(|d(t)|)$ is the essential supremum.

* By $K$ we denote the class of increasing, continuous functions $a: \mathbb{R}_+ \to \mathbb{R}_+$ with $a(0) = 0$. By $K_\infty$ we denote the class of increasing, continuous functions $a: \mathbb{R}_+ \to \mathbb{R}_+$ with $a(0) = 0$ and

$\lim_{s \to +\infty}(a(s)) = +\infty$. By $KL$ we denote the set of all continuous functions $\sigma: \mathbb{R}_+ \times \mathbb{R}_+ \to \mathbb{R}_+$ with the properties: (i) for each $t \geq 0$ the mapping $\sigma(\cdot, t)$ is of class $K$; (ii) for each $s \geq 0$, the mapping $\sigma(s, \cdot)$ is non-increasing with $\lim_{t \to +\infty}(\sigma(s, t)) = 0$.

* Let $S \subseteq \mathbb{R}^n$ be a non-empty set with $0 \in S$. We say that a function $V: S \to \mathbb{R}_+$ is positive definite if $V(x) > 0$ for all $x \in S$ with $x \neq 0$ and $V(0) = 0$. We say that a continuous function $V: S \to \mathbb{R}_+$ is radially unbounded if the following property holds: "for every $M > 0$ the set $\{x \in S : V(x) \leq M\}$ is compact". For $V \in C^1(S)$ we define $\nabla V(x) = \left( \frac{\partial V}{\partial x_1}(x), ..., \frac{\partial V}{\partial x_n}(x) \right)$.

We next recall the notions of output stability which are used in the present work. For further details the reader can consult [14].

Stability Notions for Input-Free Systems: Let $f: \mathbb{R}^n \to \mathbb{R}^n$ be a locally Lipschitz mapping and $h: \mathbb{R}^n \to \mathbb{R}^p$ be a continuous mapping with $f(0) = 0$ and $h(0) = 0$. Consider the system

$$\dot{x} = f(x), \; x \in \mathbb{R}^n \tag{1.1}$$

with output given by the equation

$$y = h(x) \tag{1.2}$$

We assume that (1.1) is forward complete, i.e., for every $x_0 \in \mathbb{R}^n$ the solution $x(t) = \phi(t, x_0)$ of the initial-value problem (1.1) with initial condition $x(0) = x_0$ exists for all $t \geq 0$. We use the notation $y(t, x_0) = h(\phi(t, x_0))$ for $t \geq 0$, $x_0 \in \mathbb{R}^n$ and $B_R = \{x \in \mathbb{R}^n : |x| < R\}$ for $R > 0$. The following stability notions are well-known (see for instance [7, 10, 28]). We say that (1.1), (1.2) is (i) *Lagrange output stable* if for every $R > 0$ the set $\{|y(t, x_0)| : x_0 \in B_R, t \geq 0\}$ is bounded, (ii) *Lyapunov output stable* if for every $\varepsilon > 0$ there exists $\delta(\varepsilon) > 0$ such that for all $x_0 \in B_{\delta(\varepsilon)}$, it holds that $|y(t, x_0)| \leq \varepsilon$ for all $t \geq 0$, (iii) *Globally Asymptotically Output Stable (GAOS)* if system (1.1), (1.2) is Lagrange and Lyapunov output stable and $\lim_{t \to +\infty}(y(t, x_0)) = 0$ for all $x_0 \in \mathbb{R}^n$, and (iv) *Uniformly Globally Asymptotically Output Stable (UGAOS)* if system (1.1), (1.2) is Lagrange and Lyapunov output stable and for every $\varepsilon, R > 0$ there exists $T(\varepsilon, R) > 0$ such that for all $x_0 \in B_R$, it holds that $|y(t, x_0)| \leq \varepsilon$ for all $t \geq T(\varepsilon, R)$.

We say that system (1.1), (1.2) is *Globally Output Attractive (GOA)* if $\lim_{t \to +\infty}(y(t, x_0)) = 0$ for all $x_0 \in \mathbb{R}^n$. We say that system (1.1), (1.2) is *practically Globally Output Attractive (p-GOA)* if there exists a constant $\tilde{\alpha} > 0$ such that $\limsup_{t \to +\infty}(|y(t, x_0)|) \leq \tilde{\alpha}$ for all $x_0 \in \mathbb{R}^n$. The constant $\tilde{\alpha} > 0$ is called the *asymptotic residual constant*.

When $h(x) = x$ then the word "output" in the above properties is omitted (e.g., Lagrange stability, Lyapunov stability, GAS, UGAS, p-UGAS, GA, p-GA).

It should be noted that (see Theorem 2.2 on page 62 in [6]) UGAOS for (1.1), (1.2) is equivalent to the existence of $\beta \in KL$ such that the following estimate holds for all $x_0 \in \mathbb{R}^n$ and $t \geq 0$:

$$|y(t,x_0)| \leq \beta(|x_0|,t) \tag{1.3}$$

We say that system (1.1), (1.2) is *practically Uniformly Globally Asymptotically Output Stable (p-UGAOS)* if there exists $\beta \in KL$ and a constant $\alpha > 0$ such that the following estimate holds for all $x_0 \in \mathbb{R}^n$ and $t \geq 0$:

$$|y(t,x_0)| \leq \beta(|x_0|,t) + \alpha \tag{1.4}$$

The constant $\alpha > 0$ is called the *residual constant*.

<u>Stability Notions for Systems with Inputs:</u> Let $D \subseteq \mathbb{R}^p$ be a closed set with $0 \in D$ and $f : \mathbb{R}^n \times D \to \mathbb{R}^n$ be a locally Lipschitz with respect to $x \in \mathbb{R}^n$ mapping with $f(0,0) = 0$. Consider the control system

$$\dot{x} = f(x,d), \, x \in \mathbb{R}^n, \, d \in D \tag{1.5}$$

We assume that (1.5) is forward complete, i.e., for every $x_0 \in \mathbb{R}^n$, $d \in L^\infty(\mathbb{R}_+;D)$ the solution $x(t) = \phi(t,x_0;d)$ of the initial-value problem (1.5) with initial condition $x(0) = x_0$ exists for all $t \geq 0$. We use the notation $y(t,x_0;d) = h(\phi(t,x_0;d))$ for $t \geq 0$, $x_0 \in \mathbb{R}^n$ and $d \in L^\infty(\mathbb{R}_+;D)$. We say that system (1.5) satisfies the *practical Uniform Bounded-Input-Bounded-State (p-UBIBS)* property if there exists a function $\bar{\gamma} \in K_\infty$ and a constant $\bar{\alpha} > 0$ such that the following estimate holds for all $x_0 \in \mathbb{R}^n$ and $d \in L^\infty(\mathbb{R}_+;D)$:

$$\sup_{t \geq 0}\left(|\phi(t,x_0;d)|\right) \leq \bar{\gamma}(|x_0|) + \bar{\gamma}(\|d\|_\infty) + \bar{\alpha} \tag{1.6}$$

The p-UBIBS property with $\bar{\alpha} = 0$ is called the UBIBS property.

We say that system (1.5), (1.2) is *Input-to-Output Stable (IOS)* if there exist $\beta \in KL$ and a non-decreasing, continuous function $\gamma : \mathbb{R}_+ \to \mathbb{R}_+$ with $\gamma(0) = 0$ such that the following estimate holds for all $x_0 \in \mathbb{R}^n$, $t \geq 0$ and $d \in L^\infty(\mathbb{R}_+;D)$:

$$|y(t,x_0;d)| \leq \beta(|x_0|,t) + \gamma(\|d\|_\infty) \tag{1.7}$$

We say that system (1.5), (1.2) is *practically Input-to-Output Stable (p-IOS)* if there exist $\beta \in KL$, a non-decreasing, continuous function $\gamma : \mathbb{R}_+ \to \mathbb{R}_+$ with $\gamma(0) = 0$ and a constant $\alpha > 0$ such that the following estimate holds for all $x_0 \in \mathbb{R}^n$, $t \geq 0$ and $d \in L^\infty(\mathbb{R}_+;D)$:

$$|y(t,x_0;d)| \leq \beta(|x_0|,t) + \gamma(\|d\|_\infty) + \alpha \tag{1.8}$$

The constant $\alpha > 0$ is called the *residual constant* while the function $\gamma$ is called the *gain function of the input* $d \in D$ *to the output* $y$.

We say that system (1.5), (1.2) satisfies the *practical Output Asymptotic Gain (p-OAG)* property if there exists a non-decreasing, continuous function $\tilde{\gamma}: \mathbb{R}_+ \to \mathbb{R}_+$ with $\tilde{\gamma}(0) = 0$ and a constant $\tilde{\alpha} > 0$ such that the following estimate holds for all $x_0 \in \mathbb{R}^n$ and for every $d \in L^\infty(\mathbb{R}_+; D)$:

$$\limsup_{t \to +\infty} \left( |y(t, x_0; d)| \right) \leq \tilde{\gamma}\left( \|d\|_\infty \right) + \tilde{\alpha} \tag{1.9}$$

The constant $\tilde{\alpha} > 0$ is called the *asymptotic residual constant* while the non-decreasing, continuous function $\tilde{\gamma}$ is called the *asymptotic gain function of the input $d \in D$ to the output $y$*. When $\tilde{\gamma} \equiv 0$ we say that system (1.5), (1.2) satisfies the *zero practical Output Asymptotic Gain property (zero p-OAG)*.

When $h(x) = x$ then the word "output" in the above properties is either replaced by the word "state" (e.g., ISS, p-ISS) or is omitted (e.g., p-AG, zero p-AG).

## 2. Main Result

In this work we study nonlinear control systems of the form

$$\dot{y} = f(y) + \sum_{i=1}^{m} g_i(y)\left( b_i u_i + \varphi_i'(y)\theta + A_i'(y)d \right)$$
$$y \in \mathbb{R}^n, d \in \mathbb{R}^q, \theta \in \mathbb{R}^p \tag{2.1}$$
$$u_1, \ldots, u_m \in \mathbb{R}, b_1, \ldots, b_m \in (0, +\infty)$$

where $f, g_1, \ldots, g_m : \mathbb{R}^n \to \mathbb{R}^n$, $\varphi_1, \ldots, \varphi_m : \mathbb{R}^n \to \mathbb{R}^p$, $A_1, \ldots, A_m : \mathbb{R}^n \to \mathbb{R}^q$ are smooth mappings with $f(0) = 0$, $\varphi_1(0) = \ldots = \varphi_m(0) = 0$, $y \in \mathbb{R}^n$ is the plant state, $u = (u_1, \ldots, u_m) \in \mathbb{R}^m$ is the control input and $\theta \in \mathbb{R}^p$, $d \in \mathbb{R}^q$, $b = (b_1, \ldots, b_m) \in (0, +\infty)^m$ are unknown disturbances. We assume next that $d \in L^\infty(\mathbb{R}_+; \mathbb{R}^q)$, $\theta \in L^\infty(\mathbb{R}_+; \mathbb{R}^p)$, $b \in L^\infty(\mathbb{R}_+; (0, +\infty)^m)$ but we assume no bounds for $\theta \in \mathbb{R}^p$, $b \in \mathbb{R}^m$ and $d \in \mathbb{R}^q$. We also assume that $\inf_{t \geq 0}(b_i(t)) > 0$ for $i = 1, \ldots, m$ but we assume no positive lower bounds for $b_i$, $i = 1, \ldots, m$.

The reader should notice the difference between $\theta$ and $d$ when both are considered to be perturbations: while $\theta \in \mathbb{R}^p$ is a vanishing perturbation (due to the fact that $\varphi_1(0) = \ldots = \varphi_m(0) = 0$), $d \in \mathbb{R}^q$ is -in general- a non-vanishing perturbation. Systems of the form (2.1) are systems that satisfy the so-called matching condition, i.e., the effect of both $\theta \in \mathbb{R}^p$ and $d \in \mathbb{R}^q$ can be cancelled by the control input $u \in \mathbb{R}$ if $\theta \in \mathbb{R}^p$, $d \in \mathbb{R}^q$ and $b \in (0, +\infty)^m$ are known.

We perform our Lyapunov-based feedback design for system (2.1) under some assumptions. Our assumptions involve an appropriate Lyapunov function.

**Assumption (A):** *There exist a constant $\sigma \in \mathbb{R}$ and smooth mappings $V, Q : \mathbb{R}^n \to \mathbb{R}_+$, $s_i : \mathbb{R}^n \to \mathbb{R}_+$, $\delta_i : \mathbb{R}^n \to \mathbb{R}$ for $i = 1,...,m$, $V, Q$ being positive definite and radially unbounded such that the following inequality holds for all $y \in \mathbb{R}^n$:*

$$\nabla V(y) f(y) \leq -Q(y) + \sigma \sum_{i=1}^{m} \delta_i(y) \nabla V(y) g_i(y) + \sum_{i=1}^{m} s_i(y) \left( \nabla V(y) g_i(y) \right)^2 \tag{2.2}$$

**Assumption (B):** *There exist a smooth mapping $\mu : \mathbb{R}^n \to (0, +\infty)$ and a constant $\Lambda \geq 0$ such that the following inequality holds for all $y \in \mathbb{R}^n$:*

$$\sum_{i=1}^{m} |\varphi_i(y)|^2 \leq \mu(y) \left( Q(y) + \Lambda \right) \tag{2.3}$$

**Remark 1:** It should be noted that Assumptions (A) and (B) are automatically valid for every system of the form (2.1) for which a smooth feedback of the form $u_i = k_i(y)$ for $i = 1,...,m$ exists that guarantees global asymptotic stability of the equilibrium point $y = 0$ when $\theta = 0$, $d = 0$ and $b_i = 1$ for $i = 1,...,m$. Notice that we do not even require that the feedback must vanish at zero, i.e., we do not require that $k_i(0) = 0$ for $i = 1,...,m$ (but $y = 0$ must be an equilibrium point of the closed-loop system (2.1) with $\theta = 0$, $d = 0$ and $b_i = 1$, $u_i = k_i(y)$ for $i = 1,...,m$, i.e., $\sum_{i=1}^{m} g_i(0) k_i(0) = 0$).

To see why Assumptions (A) and (B) are automatically valid for every system of the form (2.1) for which a smooth feedback of the form $u_i = k_i(y)$ for $i = 1,...,m$ exists, we notice that Theorem 4.5 in [1] guarantees the existence of a smooth, positive definite and radially unbounded function $V : \mathbb{R}^n \to \mathbb{R}_+$ which is a Lyapunov function for the closed-loop system (2.1) with $\theta = 0$, $d = 0$ and $b_i = 1$, $u_i = k_i(y)$ for $i = 1,...,m$, i.e., the system $\dot{y} = f(y) + \sum_{i=1}^{m} g_i(y) k_i(y)$. Furthermore, Theorem 4.5 in [1] guarantees that the inequality $\nabla V(y) f(y) + \sum_{i=1}^{m} k_i(y) \nabla V(y) g_i(y) \leq -V(y)$ holds for all $y \in \mathbb{R}^n$. It follows that (2.2) holds with $s_i(y) \equiv 0$, $\sigma = -1$, $\delta_i(y) = k_i(y)$ for $i = 1,...,m$ and $Q(y) = V(y)$ (notice that $Q$ is a smooth, positive definite and radially unbounded function). Moreover, (2.3) holds with arbitrary $\Lambda > 0$ and $\mu(y) = \Lambda + (V(y) + \Lambda)^{-1} \sum_{i=1}^{m} |\varphi_i(y)|^2$ (notice that $\mu$ is a smooth, positive function).

On the other hand, if Assumptions (A) and (B) hold with $\sigma = \Lambda = 0$ then the results that are stated below guarantee additional properties for the DADS feedback law that we construct. It should be noticed that (2.2) with $\sigma = 0$ requires that the feedback law $u_i = -s_i(y) \nabla V(y) g_i(y)$ for $i = 1,...,m$ guarantees global asymptotic stability of the equilibrium point $y = 0$ when $\theta = 0$, $d = 0$ and $b_i = 1$ for $i = 1,...,m$. This is the case of the so-called "$L_g V$"-type feedback that is frequently met in inverse optimality studies (see Chapter 4 in [2], Chapter 8 in [27] and [19, 20]). Therefore, the requirement that Assumptions (A) and (B) hold with $\sigma = \Lambda = 0$ is a restrictive requirement which cannot be guaranteed for many systems.

Consider the dynamic feedback law for $i = 1,...,m$:

$$u_i = -r_i(y,z)\nabla V(y)g_i(y) \tag{2.4}$$

$$\dot{z} = \Gamma \exp(-z)(V(y) - \varepsilon)^+ \quad , \quad z \in \mathbb{R} \tag{2.5}$$

where

$$\begin{aligned}
r_i(y,z) &:= C(\kappa + \exp(z))s_i^2(y)(\nabla V(y)g_i(y))^2 + (\kappa + \exp(z))s_i(y) \\
&\quad + C^3(\kappa + \exp(z))^3 P_i^2(y,z)(\nabla V(y)g_i(y))^2 + C(\kappa + \exp(z))^2 P_i(y,z) \\
P_i(y,z) &:= |\varphi_i(y)|^2 + |A_i(y)|^2 + (\kappa + \exp(z))^2 \mu(y) + \delta_i^2(y)
\end{aligned} \tag{2.6}$$

and $\varepsilon, \Gamma, C, \kappa > 0$ with $2C\kappa \geq 1$ are parameters of the controller (constants). We call the controller (2.4), (2.5) a Deadzone-Adapted Disturbance Suppression (DADS) controller. The controller (2.4), (2.5) combines the use of deadzone (in (2.5)) and dynamic nonlinear damping (in (2.4) and (2.6)).

The following theorem clarifies the performance characteristics that the DADS controller (2.4), (2.5) can guarantee for the closed-loop system.

**Theorem 1:** *Suppose that Assumptions (A) and (B) hold. Let $\varepsilon, \Gamma, C, \kappa > 0$ with $2C\kappa \geq 1$ be given constants. Then there exist functions $\omega \in KL$, $\zeta \in K_\infty$ and $R \in C^0(\mathbb{R}^{n+1} \times \mathbb{R}_+^2 \times (0,+\infty))$, such that for every $(y_0, z_0) \in \mathbb{R}^n \times \mathbb{R}$, $d \in L^\infty(\mathbb{R}_+; \mathbb{R}^q)$, $\theta \in L^\infty(\mathbb{R}_+; \mathbb{R}^p)$ $b \in L^\infty(\mathbb{R}_+;(0,+\infty)^m)$ with $\inf_{t \geq 0}(b_i(t)) > 0$ for $i = 1,...,m$, the unique solution of the initial-value problem (2.1), (2.4), (2.5) with initial condition $(y(0), z(0)) = (y_0, z_0)$ satisfies the following estimates for all $t \geq 0$:*

$$\limsup_{t \to +\infty}(V(y(t))) \leq \varepsilon, \tag{2.7}$$

$$V(y(t)) \leq \omega(V(y_0), t) + \zeta\left(\chi\left(\|d\|_\infty, \|\theta\|_\infty, \min_{i=1,...,m}\left(\inf_{s \geq 0}(b_i(s))\right), \exp(z_0)\right)\right), \tag{2.8}$$

$$z_0 \leq z(t) \leq \lim_{s \to +\infty}(z(s)) \leq R\left(y_0, z_0, \|d\|_\infty, \|\theta\|_\infty, \min_{i=1,...,m}\left(\inf_{s \geq 0}(b_i(s))\right)\right). \tag{2.9}$$

*where* $\chi(s_1, s_2, s_3, s_4) := \dfrac{ms_1^2 + m\sigma^2 + m\left((s_2 - \kappa - s_4)^+\right)^2 + 2C\kappa\Lambda + 2ms_3\left(\left(\dfrac{1}{s_3} - \kappa - s_4\right)^+\right)^2}{4C(\kappa + s_4)}$ *for all* $s_1, s_2, s_4 \geq 0$, $s_3 > 0$.

**Remark 2: (a)** Since $V$ is positive definite and radially unbounded, it follows from (2.7) that the zero p-OAG property from the input $d$ to the output $y$ with an asymptotic residual constant completely independent of the parameters $\theta, b$ holds. More specifically, the asymptotic residual constant $\alpha > 0$ depends only on the constant $\varepsilon > 0$.

**(b)** For constant $\theta \in \mathbb{R}^p$, $b \in (0, +\infty)^m$ inequality (2.8) shows that the p-IOS property from the disturbance $d \in L^\infty(\mathbb{R}_+; \mathbb{R}^q)$ holds for the closed-loop system (2.1), (2.4), (2.5). Indeed, the following estimate is a direct consequence of (2.8):

$$V(y(t)) \leq \omega(V(y(0)), t) + \zeta\left(\frac{m\|d\|_\infty^2}{2\kappa C}\right)$$
$$+ \zeta\left(m \frac{\sigma^2 + \left((|\theta| - \kappa)^+\right)^2 + 2\min_{i=1,\ldots,m}(b_i)\left(\left(\frac{1}{\min_{i=1,\ldots,m}(b_i)} - \kappa\right)^+\right)^2}{2\kappa C} + \Lambda\right) \quad (2.10)$$

The above estimate guarantees the p-IOS property from the disturbance $d \in L^\infty(\mathbb{R}_+; \mathbb{R}^p)$ for the closed-loop system (2.1), (2.4), (2.5). Estimate (2.7) guarantees the zero p-OAG property with an *assignable* asymptotic residual constant.

**(c)** Taking into account the above remark, we can say that in the disturbance-free case (i.e., when $d \equiv 0$), the closed-loop system (2.1), (2.4), (2.5) for constant $\theta \in \mathbb{R}^p$, $b \in (0, +\infty)^m$:
i) is Lagrange stable,
ii) satisfies the p-UGAOS property,
iii) satisfies the UGAOS property when $|\theta| \leq \kappa$, $\min_{i=1,\ldots,m}(b_i) \geq \kappa^{-1}$ and $\sigma = \Lambda = 0$,
iv) the set $\{(0, z) \in \mathbb{R}^n \times \mathbb{R}\}$ is invariant,
v) is p-GOA with an *assignable* asymptotic residual constant.

**(d)** Estimate (2.8) guarantees the following asymptotic estimate

$$\limsup_{t \to +\infty}(V(y(t))) \leq \zeta\left(\chi(l_d, l_\theta, l_b, l_z)\right) \quad (2.11)$$

where $l_d = \limsup_{t \to +\infty}(|d(t)|)$, $l_\theta = \limsup_{t \to +\infty}(|\theta(t)|)$, $l_z = \exp\left(\lim_{t \to +\infty}(z(t))\right)$ and $l_b = \min_{i=1,\ldots,m}\left(\liminf_{t \to +\infty}(b_i(t))\right)$. Thus, in the case of vanishing disturbance $\lim_{t \to +\infty}(d(t)) = 0$ and if $\sigma = \Lambda = 0$, we guarantee that at least one of the following holds:

"Either $\lim_{t \to +\infty}(y(t)) = 0$ or $\max\left(\limsup_{t \to +\infty}(|\theta(t)|), \dfrac{1}{\min_{i=1,\ldots,m}\left(\liminf_{t \to +\infty}(b_i(t))\right)}\right) > \kappa$

and $z(t) < \ln\left(\max\left(\limsup_{s \to +\infty}(|\theta(s)|), \dfrac{1}{\min_{i=1,\ldots,m}\left(\liminf_{s \to +\infty}(b_i(s))\right)}\right) - \kappa\right)$ for all $t \geq 0$"

Therefore, exact regulation of the plant state is not excluded when $\sigma = \Lambda = 0$.

**(e)** When $\sigma = \Lambda = 0$ then slight changes in the proof of Theorem 1 show that the feedback law (2.4), (2.5), (2.6) can be simplified and equation (2.6) can be replaced by the following equation:

$$\begin{aligned} r_i(y,z) &:= (\kappa + \exp(z))P_i(y,z)\left(1 + CP_i(y,z)(\nabla V(y)g_i(y))^2\right) \\ P_i(y,z) &:= s_i(y) + \dfrac{\mu(y)}{2}(\kappa + \exp(z))^2 + C(\kappa + \exp(z))\left(|A_i(y)|^2 + |\varphi_i(y)|^2\right) \end{aligned} \quad (2.12)$$

Moreover, in this case we do not need the restriction $2C\kappa \geq 1$.

When $\sigma = \Lambda = 0$ and $\varphi_i(y) \equiv 0$ for $i = 1,\ldots,m$ then slight changes in the proof of Theorem 1 show that the feedback law (2.4), (2.5), (2.6) can be simplified and equation (2.6) can be replaced by the following equation:

$$\begin{aligned} r_i(y,z) &:= (\kappa + \exp(z))P_i(y,z)\left(1 + CP_i(y,z)(\nabla V(y)g_i(y))^2\right) \\ P_i(y,z) &:= s_i(y) + C(\kappa + \exp(z))|A_i(y)|^2 \end{aligned} \quad (2.13)$$

Again, in this case we do not need the restriction $2C\kappa \geq 1$.

**(f)** Theorem 1 can be stated for less regular mappings $f, g_1, \ldots, g_m : \mathbb{R}^n \to \mathbb{R}^n$, $\varphi_1, \ldots, \varphi_m : \mathbb{R}^n \to \mathbb{R}^p$, $A_1, \ldots, A_m : \mathbb{R}^n \to \mathbb{R}^q$ but here for simplicity reasons we assume that all mappings are smooth.

## 3. Example: the uncertain double integrator

We next study the following planar system (the uncertain double integrator)

$$\begin{aligned} \dot{y}_1 &= y_2 \\ \dot{y}_2 &= \theta_1 y_1 + \theta_2 y_2 + bu + d \end{aligned} \quad (3.1)$$

where $y = (y_1, y_2) \in \mathbb{R}^2$ is the plant state, $u \in \mathbb{R}$ is the control input and $\theta = (\theta_1, \theta_2) \in \mathbb{R}^2$, $d \in \mathbb{R}$, $b \in (0, +\infty)$ are unknown disturbances. We assume next that $d \in L^\infty(\mathbb{R}_+; \mathbb{R})$, $\theta \in L^\infty(\mathbb{R}_+; \mathbb{R}^2)$, $b \in L^\infty(\mathbb{R}_+; (0, +\infty))$ but we assume no bounds for $\theta \in \mathbb{R}^2$, $b \in \mathbb{R}$ and $d \in \mathbb{R}$. We also assume that $\inf_{t \geq 0}(b(t)) > 0$ but we assume no positive lower bound for $b$.

We next design two adaptive controllers for system (3.1).

Controller (C1)

The controller (C1) is designed by assuming that $\theta = (\theta_1, \theta_2) \in \mathbb{R}^2$ and $b \in (0, +\infty)$ are constant parameters and by applying $\sigma$-modification to the controller produced by the design methodology described on pages 168-173 in the book [17] with the function

$$V(y, \hat{\theta}, \rho) = \frac{1}{2} y_1^2 + \frac{1}{2}(y_2 + cy_1)^2 + \frac{1}{2\Gamma}(\hat{\theta}_1 - \theta_1)^2 + \frac{1}{2\Gamma}(\hat{\theta}_2 - \theta_2)^2 + \frac{b}{2\Gamma}\left(\rho - \frac{1}{b}\right)^2$$

where $\Gamma, c > 0$ are constants (controller parameters) and $\hat{\theta} = (\hat{\theta}_1, \hat{\theta}_2) \in \mathbb{R}^2$, $\rho \in \mathbb{R}$ are the controller states.

Let $\bar{\sigma} > 0$ be given (arbitrary). The controller (C1) is given by the following equations

$$\begin{aligned}
u &= -\rho(y_1 + cy_2) - c\rho(y_2 + cy_1) - \rho(\hat{\theta}_1 y_1 + \hat{\theta}_2 y_2) \\
\frac{d\hat{\theta}_1}{dt} &= \Gamma y_1 (y_2 + cy_1) - \Gamma \bar{\sigma} \hat{\theta}_1 \\
\frac{d\hat{\theta}_2}{dt} &= \Gamma y_2 (y_2 + cy_1) - \Gamma \bar{\sigma} \hat{\theta}_2 \\
\dot{\rho} &= \Gamma (y_2 + cy_1)\left((1 + c^2 + \hat{\theta}_1) y_1 + (2c + \hat{\theta}_2) y_2\right) - \Gamma \bar{\sigma} \rho
\end{aligned} \tag{3.2}$$

and guarantees the following differential inequality

$$\dot{V} \leq -\min(c, \bar{\sigma}\Gamma) V(y, \hat{\theta}, \rho) + \frac{1}{2c} d^2 + \frac{\bar{\sigma}}{2}\left(\theta_1^2 + \theta_2^2 + \frac{1}{b}\right)$$

The above differential inequality guarantees the p-IOS property and the p-UBIBS property.

When $\bar{\sigma} = 0$ then the controller (C1) becomes identical with the controller produced by the design methodology described on pages 168-173 in the book [17], which guarantees the following differential inequality

$$\dot{V} \leq -cy_1^2 - \frac{c}{2}(y_2 + cy_1)^2 + \frac{1}{2c} d^2$$

which guarantees forward completeness.

The DADS controller (C2)

Setting $\rho = \kappa + \exp(z)$, employing the CLF for (3.1) $V(y) = \frac{1}{2} y_1^2 + \frac{1}{2}(y_2 + cy_1)^2$ and applying Theorem 1 and Remark 2(e), we can give the DADS controller for arbitrary $c, \varepsilon, \Gamma, C > 0$:

$$u = -\rho \left( \frac{(1-c^2)^2 + 3c^2}{2c} + \frac{\sqrt{c^2+4}+c}{c\sqrt{c^2+4}-c^2} \rho^2 + C\rho(|y|^2+1) \right)(y_2 + cy_1)$$

$$-\rho C \left( \frac{(1-c^2)^2 + 3c^2}{2c} + \frac{\sqrt{c^2+4}+c}{c\sqrt{c^2+4}-c^2} \rho^2 + C\rho(|y|^2+1) \right)^2 (y_2 + cy_1)^3 \quad (3.3)$$

$$\dot{\rho} = \Gamma \left( \frac{1}{2} y_1^2 + \frac{1}{2}(y_2 + cy_1)^2 - \varepsilon \right)^+$$

where $\rho \in (\kappa, +\infty)$ is the controller state. Indeed, the reader can verify that in this case Assumptions (A) and (B) are valid with $s(y) \equiv \frac{(1-c^2)^2 + 3c^2}{2c}$, $Q(y) = cV(y)$, $\sigma = \Lambda = 0$, $\delta(y) \equiv 0$ and $\mu(y) \equiv 2\frac{\sqrt{c^2+4}+c}{c\sqrt{c^2+4}-c^2}$.

We perform a numerical experiment for the closed-loop systems with controllers (C1), (C2). The initial conditions were

$$y_1(0) = 1, y_2(0) = 0, \rho(0) = 0.11, \hat{\theta}_1(0) = \hat{\theta}_2(0) = 0$$

and constant parameters

$$\theta_1 = \theta_2 = 1, b = 0.01, c = 0.5, \Gamma = 20, \varepsilon = 0.005, C = 1$$

The results are shown in the next pages for two particular cases: the disturbance-free case $d \equiv 0$ and the persistent disturbance case $d(t) = 2\sin(t)$.

The figures above show that no parameter drift appears for the DADS controller in the persistent disturbance case (see Fig. 10) as is the case for the controller (C1) with $\bar{\sigma} = 0$ (no leakage). Moreover, the DADS controller guarantees convergence of the plant state to a neighborhood of the origin even in the persistent disturbance case (see Fig. 3 and Fig. 8) while this is not the case for the controller (C1) with $\bar{\sigma} = 0.2$ (with leakage; see Fig. 2 and Fig. 7). In all cases the DADS controller required less control effort (see Fig. 4 and Fig. 9).

It should be noted that a comparison between Fig.1 and Fig. 8 shows that the DADS controller guarantees even in the case of a persistent disturbance almost the same performance that is guaranteed by the nominal controller (C1) with $\bar{\sigma} = 0$ (no leakage) in the disturbance-free case. Theorem 1 shows that the DADS controller achieves a practically complete disturbance suppression even in the case where $\theta = (\theta_1, \theta_2) \in \mathbb{R}^2$ and $b \in (0, +\infty)$ are time-varying.

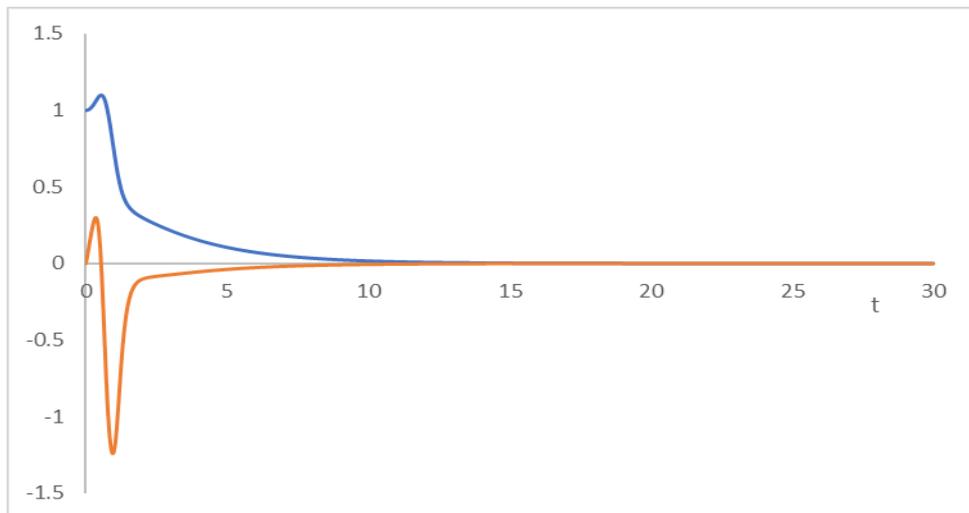

**Fig. 1:** The plant state for (C1) with $\bar{\sigma} = 0$ (no leakage) and $d \equiv 0$.
Blue line $y_1(t)$ and red line $y_2(t)$.

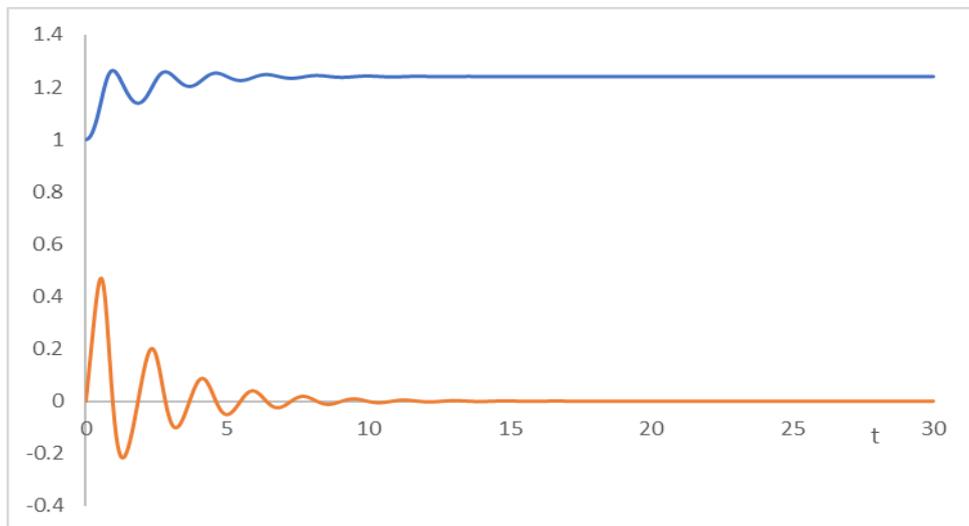

**Fig. 2:** The plant state for (C1) with $\bar{\sigma} = 0.2$ (leakage) and $d \equiv 0$.
Blue line $y_1(t)$ and red line $y_2(t)$.

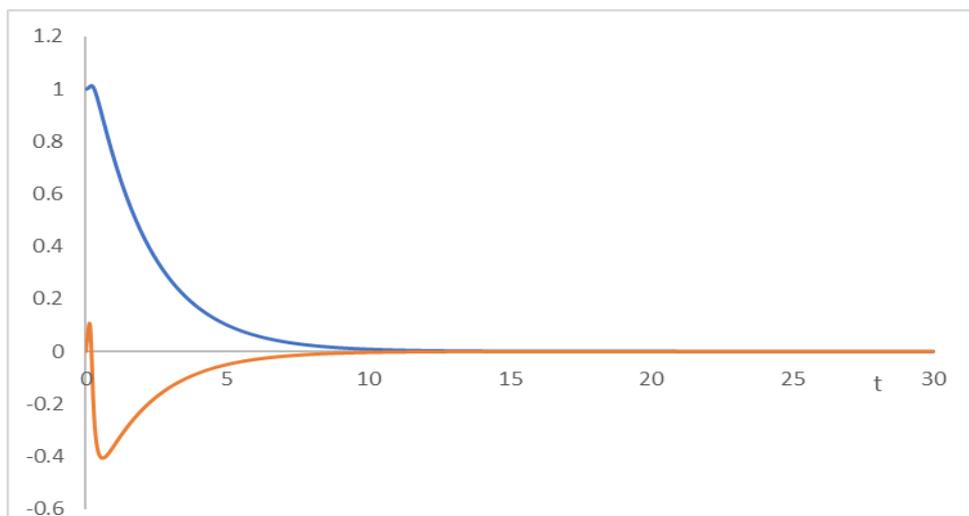

**Fig. 3:** The plant state for DADS controller (C2) with $d \equiv 0$.
Blue line $y_1(t)$ and red line $y_2(t)$.

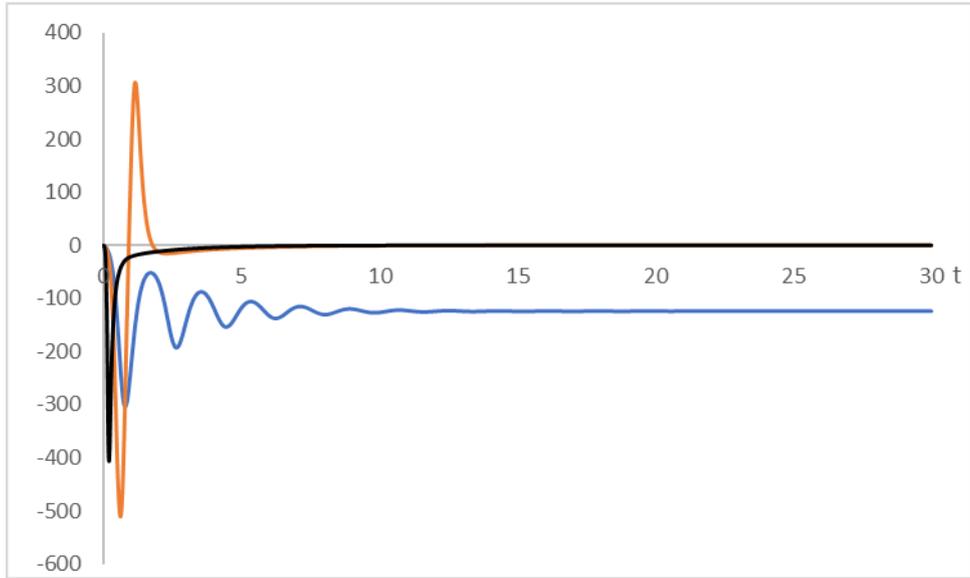

**Fig. 4:** The control input $u(t)$ with $d \equiv 0$. Red line controller (C1) with $\bar{\sigma} = 0$ (no leakage), blue line controller (C1) with $\bar{\sigma} = 0.2$ (leakage), and black line DADS controller (C2).

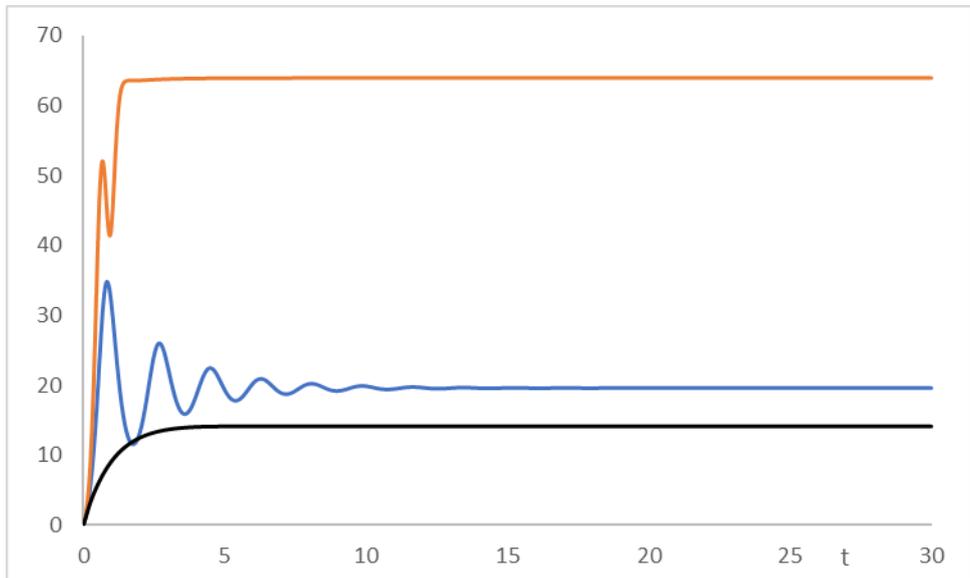

**Fig. 5:** The controller state $\rho(t)$ with $d \equiv 0$. Red line controller (C1) with $\bar{\sigma} = 0$ (no leakage), blue line controller (C1) with $\bar{\sigma} = 0.2$ (leakage), and black line DADS controller (C2).

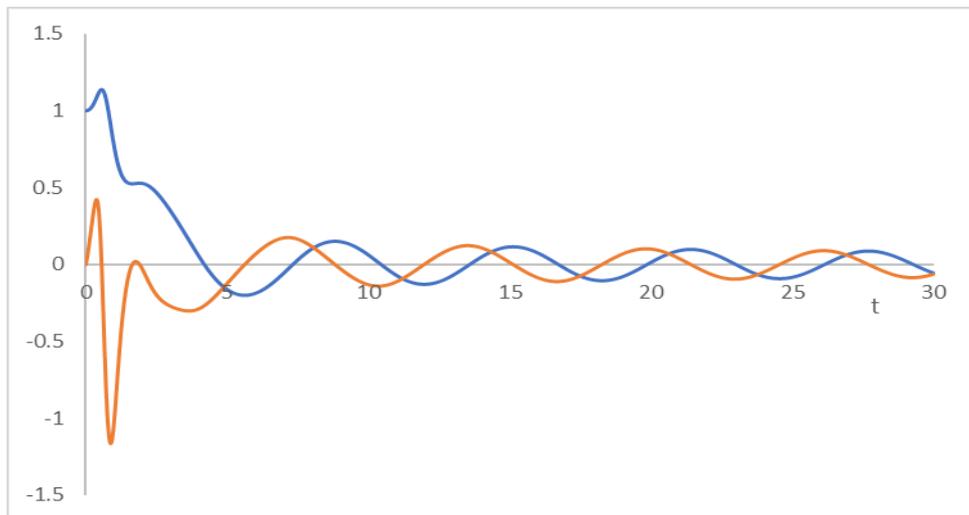

**Fig. 6:** The plant state for controller (C1) with $\bar{\sigma}=0$ (no leakage) and $d(t)=2\sin(t)$.
Blue line $y_1(t)$ and red line $y_2(t)$.

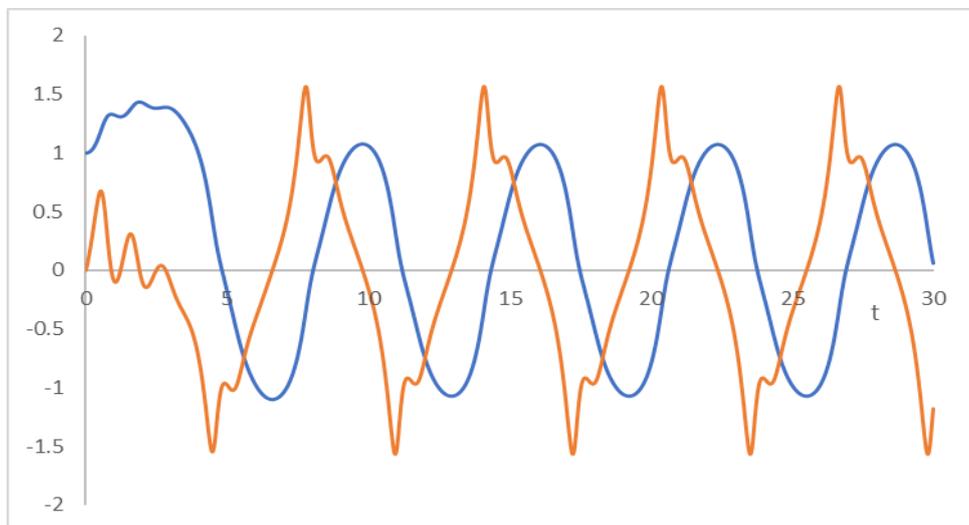

**Fig. 7:** The plant state for controller (C1) with $\bar{\sigma}=0.2$ (leakage) and $d(t)=2\sin(t)$.
Blue line $y_1(t)$ and red line $y_2(t)$.

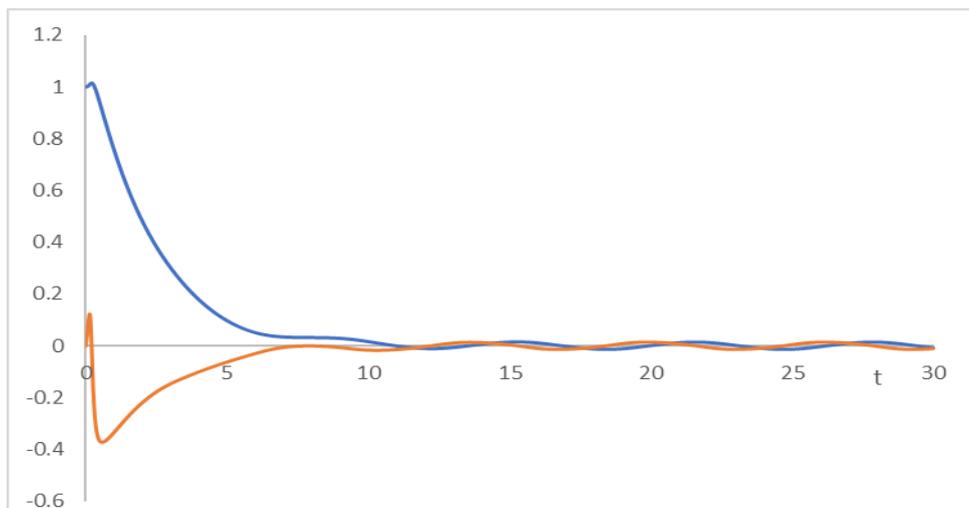

**Fig. 8:** The plant state for DADS controller (C2) with $d(t)=2\sin(t)$.
Blue line $y_1(t)$ and red line $y_2(t)$.

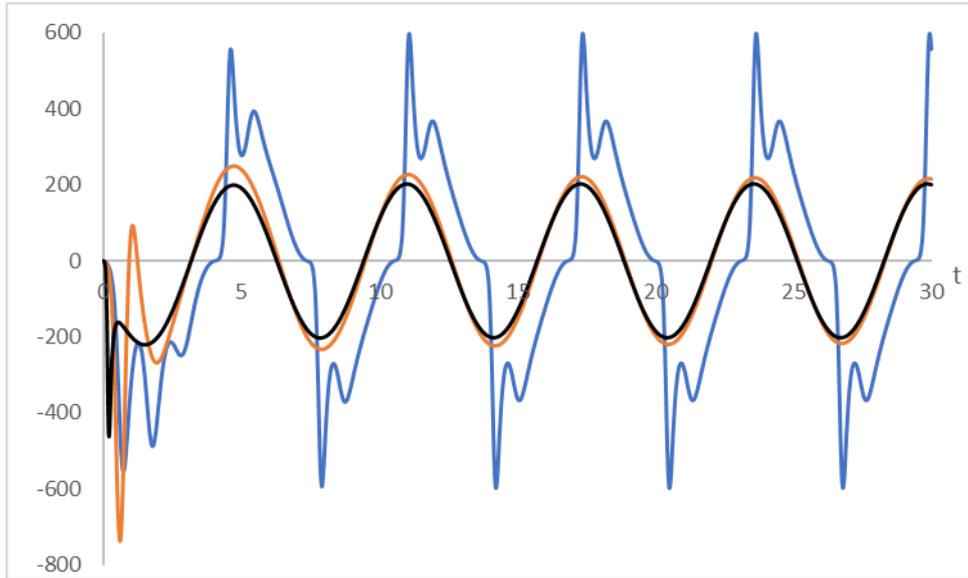

**Fig. 9:** The control input $u(t)$ with $d(t) = 2\sin(t)$. Red line controller (C1) with $\bar{\sigma} = 0$ (no leakage), blue line controller (C1) with $\bar{\sigma} = 0.2$ (leakage), and black line DADS controller (C2).

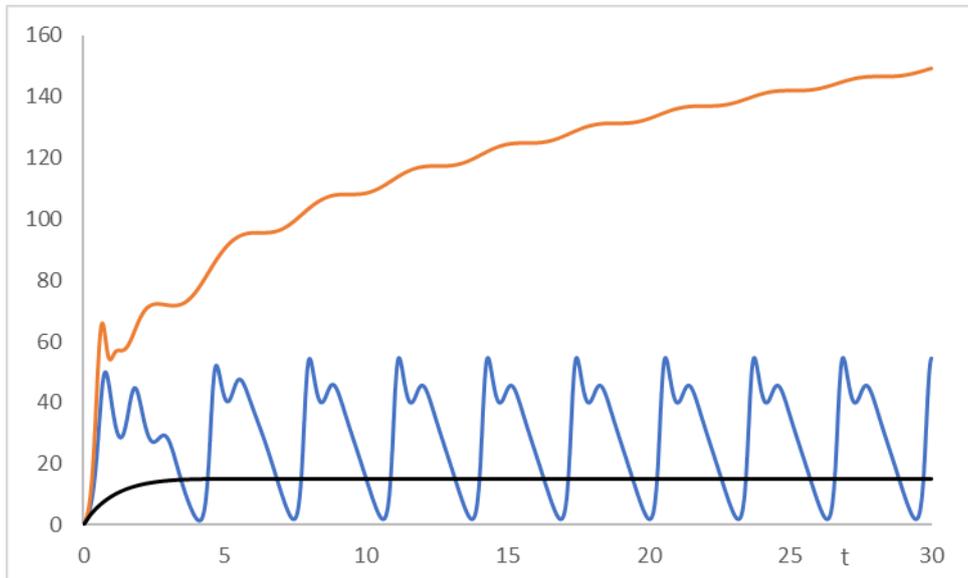

**Fig. 10:** The controller state $\rho(t)$ with $d(t) = 2\sin(t)$. Red line controller (C1) with $\bar{\sigma} = 0$ (no leakage), blue line controller (C1) with $\bar{\sigma} = 0.2$ (leakage), and black line DADS controller (C2).

## 4. Proofs of Main Results

For the proof of Theorem 1, we use the following technical lemma, proven in [15].

**Lemma 1:** *Let $\rho \in K_\infty$, $T, \varepsilon > 0$ be given. Then there exists a positive, non-increasing function $c_\varepsilon : \mathbb{R}_+ \to (0, +\infty)$ such that for every absolutely continuous function $V : [0, T) \to \mathbb{R}_+$ and for every non-increasing function $\alpha \in C^0([0, T); \mathbb{R}_+)$ for which the following differential inequality holds for $t \in [0, T)$ a.e.*

$$\dot{V}(t) \leq -\rho(V(t)) + \alpha(t) \tag{4.1}$$

*the following estimate holds for all $t_0 \in [0, T)$ and $t \in [t_0, T)$*

$$V(t) \leq \min\left(s, s\exp(-c_\varepsilon(s)(t - t_0)) + \frac{\varepsilon}{2} + \frac{\alpha(t_0)}{c_\varepsilon(s)}\right) \tag{4.2}$$

*where $s = V(0) + \rho^{-1}(\alpha(0))$.*

We next provide the proof of Theorem 1.

**Proof of Theorem 1:** Due to (2.1), (2.4) we have for all $(y, z) \in \mathbb{R}^n \times \mathbb{R}$, $d \in \mathbb{R}^q, \theta \in \mathbb{R}^p, b \in (0, +\infty)^m$:

$$\begin{aligned}\dot{V}(y, z, d, \theta, b) &= \nabla V(y) f(y) + \sum_{i=1}^m \nabla V(y) g_i(y) \varphi_i'(y) \theta \\ &+ \sum_{i=1}^m \nabla V(y) g_i(y) A_i'(y) d - \sum_{i=1}^m b_i r_i(y, z) (\nabla V(y) g_i(y))^2\end{aligned} \tag{4.3}$$

where $\dot{V}(y, z, d, \theta, b)$ is the derivative of $V(y)$ along the trajectories of the closed-loop system (2.1) with (2.4), (2.5). Using (2.2) and the inequalities for $i = 1, \ldots, m$

$$\nabla V(y) g_i(y) A_i'(y) d \leq C(\kappa + \exp(z)) |A_i(y)|^2 (\nabla V(y) g_i(y))^2 + \frac{|d|^2}{4C(\kappa + \exp(z))}$$

we get from (4.3) for all $(y, z) \in \mathbb{R}^n \times \mathbb{R}$, $d \in \mathbb{R}^q, \theta \in \mathbb{R}^p, b \in (0, +\infty)^m$:

$$\begin{aligned}\dot{V}(y, z, d, \theta, b) &\leq -Q(y) + \sigma \sum_{i=1}^m \delta_i(y) \nabla V(y) g_i(y) + \sum_{i=1}^m s_i(y) (\nabla V(y) g_i(y))^2 \\ &+ C(\kappa + \exp(z)) \sum_{i=1}^m |A_i(y)|^2 (\nabla V(y) g_i(y))^2 - \sum_{i=1}^m b_i r_i(y, z) (\nabla V(y) g_i(y))^2 \\ &+ \sum_{i=1}^m |\nabla V(y) g_i(y)| \|\varphi_i(y)\| |\theta| + \frac{m|d|^2}{4C(\kappa + \exp(z))}\end{aligned} \tag{4.4}$$

Inequality (4.4) and the inequalities

$$|\theta| \leq (|\theta| - \kappa - \exp(z))^+ + \kappa + \exp(z)$$

$$\sigma \delta_i(y) \nabla V(y) g_i(y) \leq \frac{\sigma^2}{4C(\kappa + \exp(z))} + C(\kappa + \exp(z)) \delta_i^2(y) (\nabla V(y) g_i(y))^2$$

give for all $(y,z) \in \mathbb{R}^n \times \mathbb{R}$, $d \in \mathbb{R}^q, \theta \in \mathbb{R}^p, b \in (0, +\infty)^m$:

$$\begin{aligned}
\dot{V}(y,z,d,\theta,b) \leq & -Q(y) - \sum_{i=1}^{m} b_i r_i(y,z) (\nabla V(y) g_i(y))^2 \\
& + C(\kappa + \exp(z)) \sum_{i=1}^{m} \left( |A_i(y)|^2 + \delta_i^2(y) \right) (\nabla V(y) g_i(y))^2 \\
& + \sum_{i=1}^{m} |\nabla V(y) g_i(y)| |\varphi_i(y)| (|\theta| - \kappa - \exp(z))^+ + \frac{m|d|^2 + m\sigma^2}{4C(\kappa + \exp(z))} \\
& + \sum_{i=1}^{m} s_i(y) (\nabla V(y) g_i(y))^2 + (\kappa + \exp(z)) \sum_{i=1}^{m} |\nabla V(y) g_i(y)| |\varphi_i(y)|
\end{aligned} \quad (4.5)$$

Using the inequalities

$$\begin{aligned}
& (\kappa + \exp(z)) |\nabla V(y) g_i(y)| |\varphi_i(y)| \\
& \leq \frac{\mu(y)}{2\kappa} (\kappa + \exp(z))^3 (\nabla V(y) g_i(y))^2 + \frac{\kappa}{2\mu(y)(\kappa + \exp(z))} |\varphi_i(y)|^2
\end{aligned}$$

$$\begin{aligned}
& |\nabla V(y) g_i(y)| |\varphi_i(y)| (|\theta| - \kappa - \exp(z))^+ \\
& \leq C(\kappa + \exp(z)) |\varphi_i(y)|^2 (\nabla V(y) g_i(y))^2 + \frac{\left( (|\theta| - \kappa - \exp(z))^+ \right)^2}{4C(\kappa + \exp(z))}
\end{aligned}$$

we get from (4.5) for all $(y,z) \in \mathbb{R}^n \times \mathbb{R}$, $d \in \mathbb{R}^q, \theta \in \mathbb{R}^p, b \in (0, +\infty)^m$ and for every $B > 0$ with $B \leq \min_{i=1,\ldots,m}(b_i)$:

$$\begin{aligned}
\dot{V}(y,z,d,\theta,b) \leq & -Q(y) + \frac{\kappa}{2\mu(y)(\kappa + \exp(z))} \sum_{i=1}^{m} |\varphi_i(y)|^2 \\
& + C(\kappa + \exp(z)) \sum_{i=1}^{m} \left( |A_i(y)|^2 + |\varphi_i(y)|^2 + \delta_i^2(y) \right) (\nabla V(y) g_i(y))^2 \\
& - B \sum_{i=1}^{m} r_i(y,z) (\nabla V(y) g_i(y))^2 + \frac{m|d|^2 + m\sigma^2 + m \left( (|\theta| - \kappa - \exp(z))^+ \right)^2}{4C(\kappa + \exp(z))} \\
& + \sum_{i=1}^{m} s_i(y) (\nabla V(y) g_i(y))^2 + \frac{\mu(y)}{2\kappa} (\kappa + \exp(z))^3 \sum_{i=1}^{m} (\nabla V(y) g_i(y))^2
\end{aligned} \quad (4.6)$$

Inequality (4.6) in conjunction with (2.3) and the fact that $2C\kappa \geq 1$ gives for all $(y,z) \in \mathbb{R}^n \times \mathbb{R}$, $d \in \mathbb{R}^q, \theta \in \mathbb{R}^p, b \in (0, +\infty)^m$ and for every $B > 0$ with $B \leq \min_{i=1,\ldots,m}(b_i)$:

$$\frac{\dot{V}(y,z,d,\theta,b)}{B} \leq -\frac{1}{2B}Q(y) - \sum_{i=1}^{m} r_i(y,z)(\nabla V(y)g_i(y))^2$$

$$+ C(\kappa+\exp(z))\sum_{i=1}^{m}\frac{1}{B}\Big(|A_i(y)|^2 + |\varphi_i(y)|^2 + \delta_i^2(y) + \mu(y)(\kappa+\exp(z))^2\Big)(\nabla V(y)g_i(y))^2$$

$$+ \frac{m|d|^2 + m\sigma^2 + m\big((|\theta|-\kappa-\exp(z))^+\big)^2 + 2C\kappa\Lambda}{4CB(\kappa+\exp(z))} + \sum_{i=1}^{m} s_i(y)\frac{1}{B}(\nabla V(y)g_i(y))^2$$

The fact that $\frac{1}{B} \leq \left(\frac{1}{B}-\kappa-\exp(z)\right)^+ + \kappa + \exp(z)$ in conjunction with (2.6) and the above inequality gives for all $(y,z) \in \mathbb{R}^n \times \mathbb{R}$, $d \in \mathbb{R}^q, \theta \in \mathbb{R}^p, b \in (0,+\infty)^m$ and for every $B > 0$ with $B \leq \min_{i=1,\ldots,m}(b_i)$:

$$\frac{\dot{V}(y,z,d,\theta,b)}{B} \leq -\frac{1}{2B}Q(y) - C^3(\kappa+\exp(z))^3 \sum_{i=1}^{m} P_i^2(y,z)(\nabla V(y)g_i(y))^4$$

$$+ C(\kappa+\exp(z))\sum_{i=1}^{m}\left(\frac{1}{B}-\kappa-\exp(z)\right)^+ P_i(y,z)(\nabla V(y)g_i(y))^2$$

$$- C(\kappa+\exp(z))\sum_{i=1}^{m} s_i^2(y)(\nabla V(y)g_i(y))^4 + \sum_{i=1}^{m} s_i(y)\left(\frac{1}{B}-\kappa-\exp(z)\right)^+ (\nabla V(y)g_i(y))^2 \quad (4.7)$$

$$+ \frac{m|d|^2 + m\sigma^2 + m\big((|\theta|-\kappa-\exp(z))^+\big)^2 + 2C\kappa\Lambda}{4CB(\kappa+\exp(z))}$$

Using the inequalities

$$\left(\frac{1}{B}-\kappa-\exp(z)\right)^+ s_i(y)(\nabla V(y)g_i(y))^2$$

$$\leq C(\kappa+\exp(z))s_i^2(y)(\nabla V(y)g_i(y))^4 + \frac{1}{4C(\kappa+\exp(z))}\left(\left(\frac{1}{B}-\kappa-\exp(z)\right)^+\right)^2$$

$$\left(\frac{1}{B}-\kappa-\exp(z)\right)^+ P_i(y,z)(\nabla V(y)g_i(y))^2$$

$$\leq C^2(\kappa+\exp(z))^2 P_i^2(y,z)(\nabla V(y)g_i(y))^4 + \frac{1}{4C^2(\kappa+\exp(z))^2}\left(\left(\frac{1}{B}-\kappa-\exp(z)\right)^+\right)^2$$

we obtain from (4.7) for all $(y,z) \in \mathbb{R}^n \times \mathbb{R}$, $d \in \mathbb{R}^q, \theta \in \mathbb{R}^p, b \in (0,+\infty)^m$ and for every $B > 0$ with $B \leq \min_{i=1,\ldots,m}(b_i)$:

$$\dot{V}(y,z,d,\theta,b) \leq -\frac{1}{2}Q(y) + \frac{m|d|^2 + m\sigma^2 + 2C\kappa\Lambda}{4C(\kappa + \exp(z))}$$
$$+ \frac{m\left((|\theta| - \kappa - \exp(z))^+\right)^2 + 2mB\left(\left(\frac{1}{B} - \kappa - \exp(z)\right)^+\right)^2}{4C(\kappa + \exp(z))} \quad (4.8)$$

Since $V, Q$ are positive definite and radially unbounded functions, there exists $\rho \in K_\infty$ such that

$$Q(y) \geq 2\rho(V(y)) \quad (4.9)$$

Combining (4.8) with (4.9), we obtain for all $(y,z) \in \mathbb{R}^n \times \mathbb{R}$, $d \in \mathbb{R}^q, \theta \in \mathbb{R}^p, b \in (0, +\infty)^m$ and for every $B > 0$ with $B \leq \min_{i=1,\ldots,m}(b_i)$:

$$\dot{V}(y,z,d,\theta,b) \leq -\rho(V(y)) + \frac{m|d|^2 + m\sigma^2 + 2C\kappa\Lambda}{4C(\kappa + \exp(z))}$$
$$+ \frac{m\left((|\theta| - \kappa - \exp(z))^+\right)^2 + 2mB\left(\left(\frac{1}{B} - \kappa - \exp(z)\right)^+\right)^2}{4C(\kappa + \exp(z))} \quad (4.10)$$

We next notice that there exist functions $\omega \in KL$, $\zeta \in K_\infty$ with the following property: for every $T \in (0, +\infty]$, $\tilde{\mu} \geq 0$ and for every absolutely continuous function $Y:[0,T) \to \mathbb{R}_+$ for which the following inequality holds for $t \in [0,T)$ a.e.:

$$\dot{Y}(t) \leq -\rho(Y(t)) + \tilde{\mu} \quad (4.11)$$

the following estimate holds for all $t \in [0,T)$:

$$Y(t) \leq \omega(Y(0), t) + \zeta(\tilde{\mu}) \quad (4.12)$$

The existence of functions $\omega \in KL$, $\zeta \in K_\infty$ follows from Lemma 2.14 on page 82 in [6] and by noticing that (4.11) shows that the implication

$$Y(t) \geq \rho^{-1}(2\tilde{\mu}) \Rightarrow \dot{Y}(t) \leq -\rho(Y(t))/2$$

holds for $t \in [0,T)$ a.e..

By virtue of Lemma 1 there exists a positive non-increasing function $c_\varepsilon : \mathbb{R}_+ \to (0, +\infty)$ such that for every absolutely continuous function $h:[0,T) \to \mathbb{R}_+$ and for every non-increasing function $\alpha \in C^0([0,T); \mathbb{R}_+)$ for which the differential inequality $\dot{h}(t) \leq -\rho(h(t)) + \alpha(t)$ holds for $t \in [0,T)$ a.e., the estimate

$$h(t) \leq \min\left(s, s\exp(-c_\varepsilon(s)(t-t_0)) + \frac{\varepsilon}{2} + \frac{\alpha(t_0)}{c_\varepsilon(s)}\right) \quad (4.13)$$

with $s = h(0) + \rho^{-1}(\alpha(0))$ holds for all $t_0 \in [0,T)$ and $t \in [t_0, T)$.

Moreover, define for $s \geq 0$:

$$\tilde{A}(s) := 1 + \exp(s) + \frac{2S}{\varepsilon c_\varepsilon(\bar{S})}, \tag{4.14}$$

$$G(s) := \ln\left(\tilde{A}(s) + \frac{\Gamma \bar{S}}{c_\varepsilon(\bar{S})}\right), \tag{4.15}$$

$$S := \frac{m\sigma^2 + 2ms(s+1) + 2C\kappa\Lambda}{4C}, \quad \bar{S} = s + \rho^{-1}(\kappa^{-1}S) \tag{4.16}$$

Notice that both $A$ and $G$ are non-decreasing functions.

Let arbitrary $(y_0, z_0) \in \mathbb{R}^{n+1}$, $d \in L^\infty(\mathbb{R}_+; \mathbb{R}^q)$, $\theta \in L^\infty(\mathbb{R}_+; \mathbb{R}^p)$, $b \in L^\infty(\mathbb{R}_+; (0, +\infty)^m)$ with $\inf_{t \geq 0}(b_i(t)) > 0$ for $i = 1, ..., m$ be given. Define

$$B := \min_{i=1,...,m}\left(\inf_{t \geq 0}(b_i(t))\right) \tag{4.17}$$

It follows from definition (4.17) that $B > 0$ and that $B \leq \min_{i=1,...,m}(b_i(t))$ for $t \geq 0$ a.e..

Using standard theory we conclude that there exists a unique solution $(y(t), z(t))$ of (2.1), (2.4), (2.5) with $(y(0), z(0)) = (y_0, z_0)$ defined on a maximal interval $[0, t_{max})$, where $t_{max} \in (0, +\infty]$. For every $t \in [0, t_{max})$, we have from (2.5) that $\dot{z}(t) \geq 0$. Therefore, $z(t)$ is non-decreasing on $[0, t_{max})$ and $z(t) \geq z_0$ for all $t \in [0, t_{max})$. Moreover, the mapping

$$z \mapsto \frac{m|d|^2 + m\sigma^2 + m\left(\left(|\theta| - \kappa - \exp(z)\right)^+\right)^2 + 2C\kappa\Lambda + 2mB\left(\left(\frac{1}{B} - \kappa - \exp(z)\right)^+\right)^2}{4C(\kappa + \exp(z))}$$

is non-increasing. Using (4.10), we conclude that the following inequality holds almost everywhere on $[0, t_{max})$:

$$\frac{d}{dt}(V(y(t))) \leq -\rho(V(y(t))) + \chi\left(\|d\|_\infty, \|\theta\|_\infty, B, \exp(z_0)\right) \tag{4.18}$$

Therefore, (4.11), (4.12) imply that estimate (2.8) holds for all $t \in [0, t_{max})$.

It follows from (4.10) and (4.16) with $s = \|d\|_\infty + \|\theta\|_\infty + \frac{1}{B} + V(y_0) + |z_0|$ that the following inequality holds for $t \in [0, t_{max})$ a.e.:

$$\frac{d}{dt}(V(y(t))) \leq -\rho(V(y(t))) + \frac{S}{\kappa + \exp(z(t))} \tag{4.19}$$

It follows from the fact that $z(t)$ is non-decreasing on $[0, t_{max})$, (4.19), (4.16) and (4.13) that the following estimate holds for all $t_0 \in [0, t_{max})$ and $t \in [t_0, t_{max})$:

$$V(y(t)) \leq \min\left(\bar{S}, \bar{S}\exp\left(-c_\varepsilon(\bar{S})(t-t_0)\right) + \frac{\varepsilon}{2} + \frac{S}{c_\varepsilon(\bar{S})(\kappa + \exp(z(t_0)))}\right) \tag{4.20}$$

Next, we show that $z(t)$ is bounded on $[0, t_{max})$. Since $z(t)$ is non-decreasing on $[0, t_{max})$, we have $z(t) \geq z_0$ for all $t \in [0, t_{max})$. We distinguish the following cases:

<u>*Case 1:*</u> $\exp(z(t)) \leq \tilde{A}(s)$ for all $t \in [0, t_{max})$

In this case, it follows from (4.15) that $z(t) \leq G(s)$, for all $t \in [0, t_{max})$.

<u>*Case 2:*</u> There exists $T \in (0, t_{max})$ such that $\exp(z(T)) > \tilde{A}(s)$. Since $\exp(z(0)) \leq \exp(s) < \tilde{A}(s)$, (recall (4.14) with $s = \|d\|_\infty + \|\theta\|_\infty + \frac{1}{B} + V(y_0) + |z_0|$), continuity of $z(t)$ guarantees that there exist $t_0 \in (0, T)$ such that:

$$\exp(z(t_0)) = \tilde{A}(s) \tag{4.21}$$

It follows from (4.14) that $\exp(z(t_0)) \geq \frac{2S}{\varepsilon c_\varepsilon(\bar{S})} - \kappa$ and consequently, $\frac{S}{c_\varepsilon(\bar{S})(\kappa + \exp(z(t_0)))} \leq \frac{\varepsilon}{2}$.

Therefore, we obtain from (4.20) for all $t \in [t_0, t_{max})$:

$$V(y(t)) \leq \min\left(\bar{S}, \bar{S} \exp\left(-c_\varepsilon(\bar{S})(t - t_0)\right) + \varepsilon\right) \tag{4.22}$$

Using (2.5) and inequality (4.22), we have for all $t \in [t_0, t_{max})$:

$$\frac{d}{dt}(\exp(z(t))) \leq \Gamma \bar{S} \exp\left(-c_\varepsilon(\bar{S})(t - t_0)\right) \tag{4.23}$$

Integrating (4.23) and using (4.21), we obtain for all $t \in [t_0, t_{max})$:

$$\exp(z(t)) \leq \exp(z(t_0)) + \frac{\Gamma \bar{S}}{c_\varepsilon(\bar{S})} \leq \tilde{A}(s) + \frac{\Gamma \bar{S}}{c_\varepsilon(\bar{S})} \tag{4.24}$$

Since $z(t)$ is non-decreasing, inequality (4.24) holds for all $t \in [0, t_{max})$.

It follows from (4.24) and definitions (4.14), (4.15), (4.16) that in both cases the following estimates hold for all $t \in [0, t_{max})$:

$$z_0 \leq z(t) \leq G(s) \tag{4.25}$$

$$|z(t)| \leq G(s) \tag{4.26}$$

Estimate (4.26) implies that $z(t)$ is bounded on $[0, t_{max})$.

Estimate (2.8) and the fact that $V$ is a positive definite and radially unbounded function guarantees that the component $y(t)$ of the solution $(y(t), z(t))$ of (2.1), (2.4), (2.5) with $(y(0), z(0)) = (y_0, z_0)$ is bounded on $[0, t_{max})$. Consequently, $t_{max} = +\infty$. Moreover, estimates (2.8), (4.26) hold for all $t \geq 0$. Estimates (2.8), (4.26), and the fact that $V$ is a positive definite and radially unbounded function allow us to conclude that $y \in L^\infty(\mathbb{R}_+; \mathbb{R}^n)$, $z \in L^\infty(\mathbb{R}_+; \mathbb{R})$ and $u_i \in L^\infty(\mathbb{R}_+; \mathbb{R})$ for $i = 1, ..., m$.

Since $z(t)$ is non-decreasing and bounded from above, $\lim_{\tau \to +\infty}(z(\tau))$ exists. Thus, we obtain (2.9) from (4.25) and (4.17) with

$$R\left(y_0, z_0, \|d\|_\infty, \|\delta\|_\infty, \|\theta\|_\infty, \min_{i=1,\ldots,m}\left(\inf_{t\geq 0}(b_i(t))\right)\right)$$

$$= \tilde{G}\left(\|d\|_\infty + \|\theta\|_\infty + \frac{1}{\min_{i=1,\ldots,m}\left(\inf_{t\geq 0}(b_i(t))\right)} + V(y_0) + |z_0|\right)$$

and $\tilde{G}(s) := \int_0^{s+1} G(l)\,dl$ for $s \geq 0$ (notice that $\tilde{G}: \mathbb{R}_+ \to \mathbb{R}_+$ is a non-decreasing, continuous function that satisfies $\tilde{G}(s) \geq G(s)$ for all $s \geq 0$).

We next show estimate (2.7). Estimate (4.25) and the fact that $z(t)$ is non-decreasing guarantees that the function $z(t)$ has a finite limit as $t \to +\infty$. This implies that the function $\exp(z(t))$ has a finite limit as $t \to +\infty$. Moreover, the facts that $d \in L^\infty(\mathbb{R}_+; \mathbb{R}^q)$, $\theta \in L^\infty(\mathbb{R}_+; \mathbb{R}^p)$, $y \in L^\infty(\mathbb{R}_+; \mathbb{R}^n)$, $u_i \in L^\infty(\mathbb{R}_+)$ for $i = 1, \ldots, m$, $b \in L^\infty(\mathbb{R}_+; (0, +\infty)^m)$ and (2.1) imply that $\frac{d}{dt}(V(y(t)))$ is of class $L^\infty(\mathbb{R}_+)$. It follows that the function $\frac{d}{dt}(\exp(z(t))) = \Gamma(V(y(t)) - \varepsilon)^+$ is uniformly continuous, in addition to

$$\int_0^{+\infty} \frac{d}{dt}(\exp(z(t)))\,dt = \lim_{t \to +\infty}(\exp(z(t))) - \exp(z_0) \leq \exp(G(s)) - \exp(z_0) < +\infty$$

From Barbălat's Lemma (see [10]), we have:

$$\lim_{t \to +\infty}\left(\frac{d}{dt}(\exp(z(t)))\right) = \lim_{t \to +\infty}\left(\Gamma(V(y(t)) - \varepsilon)^+\right) = 0 \qquad (4.27)$$

Therefore, estimate (2.7) holds. The proof is complete. ◁